\documentclass[12pt,leqno]{article}
\usepackage{amscd,amsmath,amsfonts,amssymb,amsthm,indentfirst}
\newcommand{\BN}{{\mathbb{N}}}
\newcommand{\BR}{{\mathbb{R}}}
\newcommand{\BC}{{\mathbb{C}}}

\newcommand{\BH}{{\mathbb{H}}}

\newcommand{\gD}{\Delta}
\newcommand{\gd}{\delta}
\newcommand{\gb}{\beta}

\newcommand{\gC}{\Gamma}
\newcommand{\gc}{\gamma}

\newcommand{\gO}{\Omega}
\newcommand{\gep}{\epsilon}

\newcommand{\ga}{\alpha}
\newcommand{\gt}{\tau}

\newtheorem{prop}{Proposition}[section]
\newtheorem{thm}[prop]{Theorem}
\newtheorem{lem}[prop]{Lemma}
\newtheorem{cor}[prop]{Corollary}
\newtheorem{conj}[prop]{Conjecture}
\theoremstyle{definition}

\newtheorem{defn}[prop]{Definition}
\newtheorem{rem}[prop]{Remark}

\title{Non-compact arithmetic manifolds have simple homotopy type 
             \\ preprint}

\author{
\sc Tsachik Gelander\\
{\small Hebrew University}\\
{\small Jerusalem 91904, Israel}\\
{\small  E-mail: tsachik@math.huji.ac.il}\\}

\begin{document} 
\maketitle

%\tableofcontents

\begin{abstract}
We formulate a conjecture that arithmetic locally symmetric manifolds
have simple homotopy type, and prove it for the non-compact case.
More precisely, we show that, for any symmetric space $S$ of non-compact type
without Euclidean de Rham factors, there are constants $\ga =\ga (S)$ and
$d=d(S)$ such that any non-compact arithmetic manifold, locally
isometric to $S$, is homotopically equivalent to a simplicial complex whose 
vertices degrees are bounded by $d$, and its number of vertices is 
bounded by $\ga$ times the Riemannian volume. 
It is very likely that such a result holds also for compact arithmetic 
manifolds.

We conclude that, for any fixed universal covering, $S$, other then the 
hyperbolic plane, there are at most $V^{CV}$ irreducible non-compact arithmetic
manifolds with volume $\leq V$, where $C=C(S)$ is a constant depending on $S$.
Since higher rank irreducible locally symmetric manifolds of finite volume
are always arithmetic, our result quantifies the number of them which are 
non-compact.
\end{abstract}

\section{Introduction}

Let $S$ be a symmetric space of non compact type without Euclidean de Rham
factors. 
An arithmetic $S$-manifold is a Riemannian manifold of the form 
$M=S/\gC$, where $\gC\leq\textrm{Isom}(S)$ is a torsion-free arithmetic lattice
in the Lie group $\textrm{Isom}(S)$ of isometries of $S$.

In this paper we prove that if $M$ is a non-compact arithmetic $S$-manifolds 
with small volume then it has a simple homotopy type. We expect that such a 
result holds also for compact arithmetic manifolds.
In fact, the compact case follows from a week version of the Lehmer's 
conjecture on the Mahler's measure of integral polynomial.
Formally, we prove the following conjecture for non-compact manifolds.

\begin{conj}
There are constants $\ga (S),d(S)$, such that any arithmetic $S$-manifold
$M=S/\gC$ with volume $\mu (M)$ is homotopically equivalent to a simplicial
complex with at most $\ga (S)\cdot\mu (M)$ vertices for which all vertices 
degrees are $\leq d(S)$.  
\end{conj}

Such information could be used in order to quantify things which are known to 
be finite, such as the number of generators or the size of a finite 
presentation for the fundamental group, as well as properties the 
homology/cohomology groups, (e.g. the Betti numbers)
 of such manifolds, in terms of their volume. 
The author's main motivation is to quantify the 
number of locally symmetric manifolds for a given universal covering.  

A classical theorem of H.C. Wang says that if $S$ is not isometric 
to one of the hyperbolic spaces $\BH^2$ or $\BH^3$, then for any $V>0$ there 
are only finitely many complete irreducible Riemannian manifolds locally 
isometric to $S$ with total volume $\leq V$, up to isometries 
(see \cite{Wa} 8.1, \cite{Bor} 8.3, and \cite{Ge} 6.5). 
We remark that Wang's result and proof do not give explicit estimates.

By Mostow's rigidity theorem, a locally symmetric manifold (with universal 
covering other then $\BH^2$) is determined by its fundamental group. Applying 
our result, we obtain a quantitive estimate on the number of non-compact
arithmetic manifolds with bounded volume.
This estimate, like the results established in \cite {BGLM} and in \cite{Ge}, 
can be considered as a step towards ``a quantitive version of Wang's theorem''.

The method used in \cite{BGLM} and in \cite{Ge} is to show that each of 
the manifolds 
under consideration has fundamental group isomorphic to the fundamental
group of a two dimensional simplicial complex with the above restrictions on
the number of vertices and on the degrees. It is not clear in general whether 
such a manifold is homotopically equivalent to a ``simple'' simplicial 
complex.
 
In \cite{BGLM} it is shown that for $n\geq 4$, the number $\rho_{\BH^n}(V)$
of complete hyperbolic $n$-manifolds with volume $\leq V$
satisfies
$$
 a_nV\log V\leq \log \rho_{\BH^n}(V)\leq b_nV\log V
$$
for some constants $a_n,b_n$ whenever $V$ is not too small.
 
In \cite{Ge} similar upper bounds are given for any rank-$1$ symmetric space
other then $\BH^2,\BH^3$, and for some examples of higher rank symmetric 
spaces (e.g. for $\BH^2\times\BH^3$ and for $\BH^2\times\BH^2\times\BH^2$).
For general $S$ other then $\BH^2,\BH^3$ and $\BH^2\times\BH^2$, theorem 3.1
in \cite{Ge} gives an upper bound for the number of compact {\it regular} 
S-manifolds.
The problem of finding upper bounds for compact non-regular manifolds is still
open. Similarly, no upper bound is known for the number of compact 
regular manifolds with universal cover $S=\BH^2\times\BH^2$.

In this note we settle the non-compact case. We show that for any $S$ other
then $\BH^2$ the number $\rho_S^{ar.,n-c}(V)$ of irreducible arithmetic 
non-compact Riemannian manifolds of type $S$ with volume $\leq V$, 
satisfies the following asymptotic estimate:
$$
 \log \rho_S^{ar.,n-c}(V)\leq C(s)V\log V.
$$
Note that, by Margulis' arithmeticity theorem, any higher rank irreducible 
locally symmetric space is arithmetic, 
while the rank-$1$ case is already dealt with in \cite{Ge}.

Our result also applies to the symmetric space $\BH^3$.
The finiteness of the number of arithmetic hyperbolic 3-manifolds with bounded
volume, as well as the finiteness of the number of arithmetic hyperbolic 
2-manifolds, was established in \cite{Bor}.

The bounds obtained here and in \cite{BGLM} and \cite{Ge} 
are independent of the symmetric space $S$, up to a normalization
of the Riemannian metric. However, in the general case, these bounds might be 
quite far from the real asymptotic behavior.
It seems, in view of \cite{Lub} , that the problem of determining the real
asymptotic behavior is closely related to the congruence subgroup problem.

\section{Statements and proofs}

Our result is:

\begin{thm}\label{thm1}
Let $S$ be a symmetric space of non compact type without Euclidean de Rham
factors. Then there are constants $\ga =\ga (S)$ and $d=d(S)$ such that if
$M=S/\gC$ is a non-compact arithmetic manifold with Riemannian volume 
$\mu (M)$, then $M$ is homotopically equivalent to a simplicial complex 
with at most $\ga\cdot\mu (M)$ vertices for which all the vertices 
degrees are bounded by $d$.
\end{thm}

A main motivation is the following application:

\begin{thm}\label{thm2}
If $S$ is not isometric to the hyperbolic plane $\BH^2$, 
then there is a constant $C=C(S)$ such that for any $V>0$, the number of 
isometric classes of non-compact irreducible arithmetic manifolds 
$M=S/\gC$ with volume $\mu (M)\leq V$ is at most $V^{C\cdot V}$.
\end{thm}

\begin{rem}
If $S$ has no $\BH^2$-factors then the irreducibility assumption in \ref{thm2}
is unnecessary. If $S$ has $\BH^2$-factors, then it is enough to require 
irreducibility in the $\BH^2$-factors. Principally, the estimation in 
\ref{thm2} follows from \ref{thm1} for all $S$-manifolds which have Mostow 
rigidity.
\end{rem}

\begin{rem}
Parts of the proof of \ref{thm1} runs along the same lines as the 
argument in \cite{BGLM} and \cite{Ge}.
We shall repeat this argument below and describe the additions and 
modifications which must be made in it in order to prove \ref{thm1}.
Some of the lemmas (e.g. \ref{positiveD}) were copied from \cite{Ge}, we bring 
them here again for the convenience of the reader.
We will also indicate how to conclude \ref{thm2} from \ref{thm1} (this is a 
word by word repetition of a short argument which already appears in 
\cite{BGLM}).
\end{rem}

We begin with an outline of the proofs of \ref{thm1} and \ref{thm2}. 
Complete details and formal definitions will be supplied later.

Fix a small enough $\gep =\gep (S)>0$,
and decompose any such $M$ into two parts $M_{\leq\gep^u}$ and 
$M_{\geq\gep^u}=\overline{M\setminus M_{\leq\gep^u}}$. In \cite{BGLM} and in
\cite{Ge} the ordinary thick-thin decomposition was used. Here we use a 
slightly different decomposition which we call {\it the unipotent thick-thin
decomposition} (see definition \ref{unip-tt} below).
We will show that $M_{\geq\gep^u}$ is a connected manifold with boundary, 
homotopically equivalent to $M$. We then show how to modify 
$M_{\geq\gep^u}$ to the subset $M_{\geq^u}=M\setminus (M_{\leq\gep^u})_\gep$, 
which is homotopically the same as $M_{\geq\gep^u}$, but on which we have
some control on the curvature of the boundary. 
We will show that the injectivity radius at any point 
in the modified unipotent thick part $M_{\geq^u}$ is at least $\gep /m$ for 
some fixed constant $m=m(S)\in\BN$. We will also show that for any point on 
the boundary of the modified unipotent thin part, there is a unit tangent 
vector with respect to which 
the directional derivative of the distance function from the unipotent thin
part is homogeneously bounded. As shown in detail in 
\cite{Ge}, this geometric information enables us to construct a simplicial 
complex within the modified unipotent
thick part which is homotopically equivalent to it, in such a way so that the
volume information $\mu (M)$ translates to the desired combinatorial 
conditions on this simplicial complex. 
The construction of such a simplicial complex is done by finding a 
good cover (in the sense of \cite{BT}) for the modified unipotent thick part 
$M_{\geq^u}$ which consist of at most $\ga\cdot\mu (M)$ sets, each of
them intersects at most $d$ of the others. The nerve of this cover constitutes 
a simplicial complex which satisfied the conditions stated in \ref{thm1}.

Mostow's rigidity theorem implies that non-isometric locally symmetric 
manifolds have non-isomorphic fundamental groups. The fundamental group of a 
simplicial complex is the one of its 2-skeleton. 
A rough estimation on the number of 2-dimensional simplicial complexes 
with bounded degrees and with bounded number of vertices, 
yields \ref{thm2}.  
\\ \\
 
We start with a variant of the Margulis lemma.
Let $G$ denote the Lie group of isometries of $S$,
$$
 G=\textrm{Isom}(S).
$$
$G$ is center-free, semi-simple without compact factors, and with 
finitely many connected components. Hence there is an algebraic structure, 
with respect to which $G$ is a real algebraic group.
For $g\in G$ we denote by $d_g:S\to \BR$ the displacement function
$$
 d_g(x)=d(x,g\cdot x).
$$
In what follows $\mu$ denotes a fixed Haar measure on $G$, as well as the 
Riemannian measure on $S$.

\begin{lem}\label{uni-thin}
There are constants $\gep =\gep (S)>0$ and $m=m(S)\in\BN$, such that 
if $\gC\leq G$ is a non-uniform torsion-free arithmetic lattice, then 
for any $x\in S$, the group of real points of the Zariski closure
$\overline{\gC_\gep (x)}^z$ of the group 
$$
 \gC_\gep (x)=\langle\gc\in\gC:d_\gc (x)\leq\gep\rangle
$$
has at most $m$ connected components, and its identity component is a 
unipotent group.
\end{lem}

\begin{proof}
Let $\gO_1\subset G$ be a Zassenhaus neighborhood (see \cite{Rag} definition
8.22 and theorem 8.16), 
and let $\gO_2\subset G$ be an identity neighborhood for which the 
intersection $\gO_2\cap\gD$ consists of unipotent elements only, 
for any non-uniform arithmetic lattice $\gD\leq G$ (see \cite{Mar1} 4.21 
page 322).
Let $\gO\subset G$ be a relatively compact symmetric identity neighborhood 
which satisfies 
$\gO^2\subset\gO_1\cap\gO_2$.

Fix 
$$
 m> \inf_{h\in G}\frac{\mu (\{g\in G:d_g(x)\leq 1\}\cdot h\gO h^{-1} )}
 {\mu (\gO )}, 
$$ 
and
$$
 \gep =\frac{1}{m}.
$$ 
As $G$ is unimodular, $m$ is independent of $x$. 
Replacing $\gO$ by its conjugation $h\gO h^{-1}$, if needed, we can assume that
$$
 m>\frac{\mu (\{g\in G:d_g(x)\leq 1\}\cdot\gO )}{\mu (\gO )}.
$$ 
Let
$$
 \gC_\gO=\langle\gC\cap\gO^2\rangle.
$$
Then 
$$
 [\gC_\gep (x): \gC_\gep (x)\cap\gC_\gO ]\leq m.
$$
To see this , assume for a moment that this index was $\geq m+1$. Then we 
could find $m+1$ representative $\gc_1,\gc_2,\ldots,\gc_{m+1}\in\gC_\gep (x)$
for different cosets of $\gC_\gep (x)\cap\gC_\gO$ in the ball of radius $m$
in $\gC_\gep (x)$ according to the word metric with respect to the generating
set $\{\gc\in\gC_\gep (x): d_\gc (x)\leq\gep\}$. As they belong to different
cosets, $\gc_i\gO\cap\gc_j\gO =\emptyset$ for any $1\leq i<j\leq m+1$.
Since $d_{\gc_i}(x)\leq m\cdot\gep =1$ these $\gc_i$'s
are all in $\{ g\in G:d_g(x)\leq 1\}$, this contradicts the assumption 
$m\cdot\mu (\gO )>\mu (\{ g\in G:d_g(x)\leq 1\}\cdot\gO )$.

It follows from the Kazhdan-Margulis theorem (see \cite{Rag} theorem 8.16)
that $\gC_\gO$ is contained in a 
connected nilpotent Lie subgroup of $G$, and therefore, by Lie's theorem,
$\gC_\gO$ is triangulizable over $\BC$. As $\gC_\gO$ is generated by unipotent
elements, it follows that $\gC_\gO$ is a group of unipotent elements.
Thus the Zariski closure $\overline{\gC_\gO}^z$ is unipotent algebraic group, 
and the group of its real points $(\overline{\gC_\gO}^z)_\BR$ is connected in 
the Hausdorff topology. Similarly, its subgroup
$(\overline{\gC_\gep (x)\cap\gC_\gO}^z)_\BR$ is connected.
Clearly $(\overline{\gC_\gep (x)\cap\gC_\gO}^z)_\BR$ is the identity 
connected component of $(\overline{\gC_\gep (x)}^z)_\BR$, 
and its index is at most $m$.

\end{proof}

\begin{rem}
Although $\gep$ could be taken to be $1/m$, we use different letters for them
because they play different roles.
In the sequel we will assume that the above lemma is satisfied with $\gep$
replaced by $10\gep$.
\end{rem}

For $g\in G$ we denote by $\{ d_g\leq\gep\}$ the sub-level set
$$
 \{ d_g\leq\gep\} =\{ x\in S:d_g(x)\leq\gep\}.
$$
More generally, for any function $\phi :L\to \BR$ defined on an abstract set 
$L$, and for any $\gt\in\BR$, we let $\{\phi\leq\gt\}$ denote
$$
 \{\phi\leq\gt\} =\{ x\in L:\phi (x)\leq\gt\}.
$$
We denote by $\gC^u$ the set of unipotent elements in $\gC$,
$$
 \gC^u=\{ \gc\in\gC :\gc \textrm{ is unipotent}\}.
$$

\begin{defn}[The unipotent thick-thin decomposition for $M=S/\gC$]
\label{unip-tt}
In the universal covering $S=\tilde M$ we take:
$$
 \tilde M_{\leq\gep^u}=\cup_{\gc\in\gC^u\setminus\{ 1\}}\{ d_\gc\leq\gep\},
$$
and
$$
 \tilde M_{\geq\gep^u}=\overline{\tilde M \setminus \tilde M_{\leq\gep^u}}.
$$
We define the unipotent thin part $M_{\leq\gep^u}\subset M$ and the unipotent 
thick part $M_{\geq\gep^u}\subset M$ to be the images of 
$\tilde M_{\leq\gep^u}$ and $\tilde M_{\geq\gep^u}$ under the universal 
covering map.\\
\end{defn}

It follows from lemma \ref{uni-thin} that if $\gc$
is an element of a non-uniform arithmetic lattice of $G$ and 
$\inf d_\gc < \gep$, then $\gc^j$ is unipotent for some $j\leq m$.
Since $d_{\gc^j}(x)\leq j\cdot d_\gc (x)$, we obtain:

\begin{cor}\label{InjRad}
The injectivity radius at any point of $M_{\geq\gep^u}$ is at least $\gep /m$.
\end{cor}

We whish to replace $M=S/\gC$ by its unipotent thick part $M_{\geq\gep^u}$ 
(where the injectivity radius is large), and then to modify $M_{\geq\gep^u}$
(to a subset with ``smother'' boundary) without changing the homotopy type. 
We shall use the following lemma twice:

\begin{lem}\label{deformation-retract}
Let $\mathcal{F} =\{\phi_i\}_{i\in I}$ be a family of non-negative 
continuous functions on $S$, such that for any $x\in S$ and any $\gt >0$ 
the set $\Psi_{x,\gt } =\{\phi\in\mathcal{F} : \phi (x)\leq\gt\}$ is 
finite, and such that each $\phi\in\mathcal{F}$ is $C^1$ on the set 
$\overline{\{\phi >0\}}$. (The gradient $\nabla\phi$ of $\phi$ on the 
boundary $\partial\{\phi >0\}$ is defined as the tangent vector with 
length and direction equal to the value of the maximal directional 
derivative and the direction at which it occurs.) 
Let $L=\cap_{\phi\in\mathcal{F}}\overline{\{\phi >0\}}$.
Assume that there is a continuous function $\gb :\BR^{\geq 0}\to\BR^{>0}$ 
(in case all $\phi\in\mathcal{F}$ are strictly positive we allow $\gb$ to
be defined only on $\BR^{>0}$)
such that for any subset $\Psi\subset\mathcal{F}$ for which 
$\cap_{\phi\in\Psi}\{\phi\leq 3\gep\}\neq\emptyset$, and for any $x\in L$, 
there is a unit tangent vector $\hat n(x,\Psi )\in T_x(S)$ such that  
$$
 \hat n(x,\Psi )\cdot \nabla \phi (x)\geq\gb\big( \phi (x)\big)
$$
for any $\phi\in\Psi$. 
Then there is a deformation retract from $L$ to the ``$\mathcal{F}$-thick 
part'' $L_{\geq\gep^\mathcal{F}}=\cap_{i\in I}\{\phi\geq\gep\}$. (In particular
it follows that $L_{\geq\gep^\mathcal{F}}\neq\emptyset$.)

If $L$ is $S=\tilde M$ or $\tilde M_{\geq\gep^u}$, and the family 
$\mathcal{F}$ is $\gC$-invariant, in the sense that the function 
$\gc\cdot\phi (x)=\phi(\gc^{-1}\cdot x)$ belongs to $\mathcal{F}$ for any
$\phi\in\mathcal{F},~\gc\in\gC$,
then there exist such a deformation retract
from $M=S/\gC$ (resp. from $M_{\geq\gep^u}=\tilde M_{\geq\gep^u}/\gC$) to the
image of $L_{\geq\gep^\mathcal{F}}$ under the universal covering map.
\end{lem}

\begin{proof}
We will define an appropriate continuous vector field on $L$. 
The desired deformation retract will be the flow along this vector field.

For any non-empty subset $\Psi\subset\mathcal{F}$ for which 
$\cap_{\phi\in\Psi}\{\phi\leq 3\gep\}\neq\emptyset$ and any $x\in L$, let 
$\hat f(x,\Psi )\in T_x(S)$ be a unit tangent vector which maximizes the 
expression
$$
 \min\{\hat f\cdot \nabla\phi (x):\phi\in\Psi\}.
$$
It follows from the strictly convexity of the Euclidean unit disk that
$\hat f(x,\Psi )$ is uniquely determined, and consequently, that for a fixed
$\Psi$, the vector field $\hat f(x,\Psi )$ is continuous.
Additionally 
$$
 \hat f(x,\Psi )\cdot \nabla \phi (x)\geq
 \hat n(x,\Psi )\cdot \nabla \phi (x)\geq\gb\big( \phi (x)\big)
$$
for any $\phi\in\Psi$.

Let $\gd (x)$ denote 
$$
 \gd (x)=\min_{\phi\in \mathcal{F}}\phi (x).
$$
The desired vector field is defined as follows:
\begin{eqnarray*}
 &&\!\!\!\!\!\!\! \overrightarrow V(x)=
 \sqrt{2\big(\gep -\gd(x)\big)\vee 0}\cdot\\
&\phantom{\le}&
 \cdot\sum_\Psi
 \frac{\big( 3\gep -\max_{\phi\in\Psi}\phi (x) \big)\vee 0}{\gep}\cdot
 \frac{\big(\min_{\phi\notin\Psi}\phi (x)-\gep \big)\vee 0}{\gep}\cdot
 \frac{1}{\gb \big(\gd (x)\big)}
 \hat f(x,\Psi ),
\end{eqnarray*}
where the sum is taken over all non-empty finite subsets 
$\Psi \subset{\cal F}$.

Clearly all the coefficients are continuous, and $\overrightarrow V\equiv 0$
on the $\mathcal{F}$-thick part
$$
 \{\gd\geq\gep\}=\cap_{\phi\in \mathcal{F}}\{\phi\geq\gep\}.
$$
The term $\sqrt{2\big(\gep -\gd(x)\big)\vee 0}$ takes care of the continuity
on the boundary $\{\gd =\gep\} =\partial\{\gd \leq\gep\}$ of the 
$\mathcal{F}$-thin part. 
The terms $\frac{\big( 3\gep -\max_{\phi\in\Psi}\phi (x) \big)\vee 0}{\gep}$ 
guarantee that all the non-zero summands correspond to sets which are 
contained in the finite set $\Psi_{x,3\gep}$. In particular the summation is 
finite for any $x\in L$, and $\hat f(x,\Psi )$ is defined for any 
non-zero summand.
The terms $\frac{\big(\min_{\phi\notin\Psi}\phi (x)-\gep \big)\vee 0}{\gep}$ 
guarantee that all the non-zero summands correspond to $\Psi$'s which contains
$\Psi_{x,\gep}$.   

If $\phi\in\mathcal{F}$ satisfies $\phi (x)=\gd (x)$ then
$$
 \nabla \phi (x)\cdot \hat f(x,\Psi )\geq \gb \big(\gd (x)\big),
$$
for any $\Psi$ which contains $\phi$. Thus, if $\gd (x)<\gep$ then
$$
 \nabla \phi (x)\cdot\overrightarrow V(x)\geq\sqrt{2\big(\gep -\gd(x)\big)}.
$$
To see this, we need only to look at the summand which corresponds to 
$\Psi =\Psi_{x,2\gep}$.

It follows that if $x(t)$ is an integral curve of $\overrightarrow V$ with
$\gd\big( x(0)\big) <\gep$, then
$$
 \frac{d}{dt}\Big(\gd\big( x(t)\big)\Big)\geq 
 \sqrt{2\Big(\gep -\gd\big( x(t)\big)\Big)}.
$$
(To be more precise, since $\gd\big( x(t)\big)$ is not differentiable, 
we should write \\
$\liminf_{\gt\to 0}\frac{\gd\big( x(t+\gt)\big)-\gd\big( x(t)\big)}{\gt}$ 
instead of $\frac{d}{dt}\Big(\gd\big( x(t)\big)\Big)$ in the estimation 
above.)
Thus, for $t=\sqrt{2\Big(\gep -\gd\big( x(0)\big)\Big)}$ we have
$\gd\big( x(t)\big) =\gep$. 

Since $x\in L$ belongs to $\partial L$ iff $\gd (x)=0$, and hence, 
the vector field 
$\overrightarrow V$ pointing everywhere towards the interior $\textrm{int}(L)$,
it follows from Peano existence theorem of solution for ordinary differential
equations, that for any $x\in L$ there is an integral curve $x(t)$ of 
$\overrightarrow V$, defined for all $t\geq 0$ with $x(0)=x$ and with 
$x(t)\in\textrm{int} L$ for $t>0$.

Concluding the above discussion, we get that the flow along 
$\overrightarrow V$ for
$\sqrt{2\gep}$ time units defines a deformation retract from $L$ to 
$=L_{\geq\gep^\mathcal{F}}=\{\gd\geq\gep\}.$

If $L$ is $S=\tilde M$ or $\tilde M_{\geq\gep^u}$ and $\mathcal{F}$ is 
$\gC$-invariant, then $\overrightarrow V(x)$, as it is defined above, is also 
$\gC$-invariant. Hence it induces a vector field on $L/\gC$, 
and a deformation retract from $L/\gC$ to $L_{\geq\gep^{\mathcal{F}}}/\gC$.
\end{proof}

We will apply lemma \ref{deformation-retract} to prove propositions
\ref{deformation-retract-1} and \ref{deformation-retract-2}.

\begin{prop}\label{deformation-retract-1}
There is a deformation retract from $M$ to
$M_{\geq\gep^u}$.
\end{prop}

\begin{proof}
We will show that the conditions of lemma \ref{deformation-retract} are 
satisfied with $L=S$ and $\mathcal{F}=\{d_\gc\}_{\gc\in\gC^u\setminus\{ 1\}}$. 
The finiteness of the sets $\Psi_{x,\gt }=\{\gc\in\gC^u\setminus\{ 1\}:
d_\gc (x)\leq\gt\}$ follows from the compactness of $\{ g\in G:d_g(x)\leq\gt\}$
together with the discreteness of $\gC$. All the functions 
$\{ d_\gc\}_{\gc\in\gC^u\setminus\{ 1\}}$ are strictly positive. 
We shell find a continuous function $\gb :\BR^{>0}\to\BR^{>0}$, and an 
appropriate direction $\hat n(x)\in T_x(S)$ for any $x\in S$.

Let $\gc_1,\gc_2,\ldots,\gc_k$ be unipotent elements of $\gC$ for which 
$\cap_{i=1}^k \{d_{\gc_i} <3\gep\}\neq \emptyset$.
By lemma \ref{uni-thin} the Zariski closure $\overline\gD^z$ of
the group $\gD =\langle \gc_1,\gc_2,\ldots,\gc_k \rangle$ 
has a unipotent identity component. Since $\gc_i$ is unipotent,
it is contained in the Zariski closure of the cyclic group generated by any 
power of it. 
As $\gc_i^j$ belongs to the identity component $(\overline{\gD}^z)^0$ for some
$j\leq m$, $\overline{\gD}^z$ is Zariski connected and hence the group of its
real points $(\overline{\gD}^z)_\BR$ is a connected unipotent group which 
contains $\gc_1,\gc_2,\ldots,\gc_k$.

Let $N\leq G$ be a maximal 
connected unipotent subgroup which contains $\gc_1,\gc_2,\ldots ,\gc_k$.
Let $W\leq S(\infty )$ be the Weyl chamber of the Tits boundary of $S$ 
which corresponds to $N$. Fix arbitrarily $x\in S$ and let $c(t)=c_x(t)$ be the
geodesic line with $c(0)=x$ for which $c(-\infty )$ is the center of $W$.

It is easy to see that the stabilizer of $c(\infty )$ intersects $N$ 
trivially. Therefore 
$$
 \frac{d}{dt}|_{t=0}\Big( d_g\big( c(t)\big)\Big)>0
$$ 
for any $g\in N\setminus\{ 1\}$. In addition, the continuous function
$$
 h(g)=\dot c(0)\cdot\nabla d_g(x)
$$ 
attains a minimum on the compact set
$$
 \{g\in N:d_g(x)=\gt\}.
$$
We define $\gb (\gt )$ to be this minimum. 
Clearly, $\gb$ is a 
continuous positive function independent of $\gc_1,\gc_2,\ldots,\gc_k$.
The conditions of lemma \ref{deformation-retract} are satisfied with the 
tangent vector $\hat n(x,\{\gc_1,\gc_2,\ldots ,\gc_k\} )=\dot c_x(0)$.
\end{proof}

We require few more notations:\\
For a subset $A$ of $S$ or of $M$, we denote by $D_A(x)$ the distance 
function
$$
 D_A(x)=d(x,A),
$$
and by $(A)_\gt$ its $\gt$-neighborhood
$$
 (A)_\gt =\{x:D_A(x)<\gt\} .
$$

We wish to replace the unipotent thick-thin decomposition 
$$
 M=M_{\leq\gep^u}\cup M_{\geq\gep^u}
$$ 
by the ``more smooth'' decomposition
$$
 M=M_{\leq^u}\cup M_{\geq^u}
$$
where $M_{\leq^u}=(M_{\leq\gep^u})_\gep$ and $M_{\geq^u}=
M\setminus (M_{\leq\gep^u})_\gep$
for which we have a control on the curvature of the boundary (see 
\ref{smallC}). In order to do this we need the following:

\begin{prop}\label{deformation-retract-2}
There is a deformation retract from $M_{\geq\gep^u}=M\setminus 
M_{\leq\gep^u}$ to $M_{\geq^u}=M\setminus (M_{\leq\gep^u})_\gep$.
\end{prop}

\begin{proof}
We just have to show that the conditions of lemma \ref{deformation-retract}
are satisfied with $L=\tilde M_{\geq\gep^u}$ and 
${\cal F}=\{ D_{\{ d_\gc\leq\gep\}}\}_{\gc\in\gC^u\setminus\{ 1\}}$.
Again, the finiteness of the sets $\Psi_{x,\gt}=\{\gc\in\gC^u\setminus\{ 1\}:
D_{\{d_\gc\leq\gep\}}(x)\leq\gt\}$ follows from the discreteness of $\gC$
together with the compactness of $\{g\in G:D_{\{d_g\leq\gep\}}(x)\leq\gt\}$. 
We shell define the direction $\hat n(x,\Psi )\in T_x(S)$ analogously 
to the way it is done in the proof of \ref{deformation-retract-1}, 
and shell show that the conditions of lemma \ref{deformation-retract} are 
satisfied with the constant function $\frac{\gb (\gep )}{2}$, where $\gb$ 
is the function defined in the proof of \ref{deformation-retract-1}.

Let $\gc_1,\gc_2,\ldots,\gc_k\in\gC$ be unipotent elements for which 
$\cap\{ D_{\{ d_{\gc_i}\leq\gep\}}\leq 3\gep\}\neq\emptyset$. Since 
$$ 
 (\{ d_{\gc_i}\leq\gep\} )_{3\gep}\subset\{d_{\gc_i}\leq 7\gep\}
$$
it follows that $\gc_1,\gc_2,\ldots,\gc_k$ are contained in a connected 
unipotent group. As in the proof of \ref{deformation-retract-1}, let $N$ be a 
maximal connected unipotent group which contains $\gc_1,\gc_2,\ldots,\gc_k$, 
let $W\subset S(\infty )$ be the Weyl chamber which corresponds to $N$ in the 
Tits boundary, and for a fixed $x\in\tilde M_{\geq\gep^u}$ let 
$c(t)=c_x(t)$ be the geodesic line with $c(0)=x$ which comes from the center 
$c(-\infty)$ of $W$. Since $d_{\gc_i}(x)>\gep$ and since
$d\big( c(t),\gc_i\cdot c(t)\big)$ tends to $0$ as $t\to -\infty$, we have
$d_{\gc_i}\big( c(t_0)\big)=\gep$ for some negative $t_0$. 

Since the set $\{ d_{\gc_i}\leq\gep\}$ is convex, the distance function
$D_{\{d_{\gc_i}\leq\gep\}}$ is convex. Therefore the function
$D_{\{d_{\gc_i}\leq\gep\}}\big( c(t)\big)$ has non-decreasing derivative,
Thus, taking $\hat n(x,\{\gc_1,\gc_2\ldots ,\gc_k\} )=\dot c(0)$, we have
\begin{eqnarray*}
  \hat n(x,\{\gc_1,\gc_2\ldots ,\gc_k\} )\cdot\nabla D_{\{d_{\gc_i}\leq\gep\}}
 (x)&=&\\
&=& \dot c(0)\cdot\nabla D_{\{d_{\gc_i}\leq\gep\}}(x)\\ 
&=& \frac{d}{dt}|_{t=0}D_{\{d_{\gc_i}\leq\gep\}}\big( c(t)\big)\\
& \geq&  
 \frac{d}{dt}|_{t=t_0}D_{\{d_{\gc_i}\leq\gep\}}\big( c(t)\big)\\
& =&
 \dot c(t_0)\cdot\nabla D_{\{d_{\gc_i}\leq\gep\}}\big( c(t_0)\big)\\
& =& 
 \dot c(t_0)\cdot\frac{\nabla d_{\gc_i}\big( c(t_0)\big)}
 {\|\nabla d_{\gc_i}\big( c(t_0)\big)\|}\\
&\geq &
 \dot c(t_0)\cdot\frac{\nabla d_{\gc_i}\big( c(t_0)\big)}{2}
 \geq\frac{\gb(\gep)}{2},
\end{eqnarray*}
where $\gb$ is the function defined in the proof of 
\ref{deformation-retract-1}.
Since $\|\nabla  D_{\{d_{\gc_i}\leq\gep\}}\| =1$ everywhere outside 
$\{d_{\gc_i}<\gep\}$, 
and $\|\nabla d_{\gc_i}\|\leq 2$ as $d_{\gc_i}$ is 2-Lipschitz.
\end{proof}

Set $b=\frac{2}{\gb(\gep )}$. For any $x\in S$ we denote
$$
 \Psi_{x,\gt}=\{\gc\in\gC^u\setminus\{ 1\} :D_{\{\gc\leq\gep\}}(x)\leq 
 \gt\},
$$
and we chose the direction $\overrightarrow F(x)\in T_x(S)$ 
which maximizes the expression 
$$
 \min_{\gc\in\Psi_{x,2\gep}} \overrightarrow F\cdot\nabla
 D_{\{\gc\leq\gep\}}(x).
$$ 
Then as it is shown in the proof of \ref{deformation-retract-2},
$$
 \overrightarrow F(x)\cdot\nabla D_{\{\gc\leq\gep\}}(x)>1/b,
$$ 
whenever $\gc\in\Psi_{x,2\gep}$ and $x\in\tilde M_{\geq\gep^u}=
\tilde M\setminus(\tilde M_{\leq\gep^u})_\gep$. 
This implies:

\begin{prop}\label{closeN}{}
For any $\gep >t>0$ the Hausdorff distance 
$$
 \textrm{Hd}\big(\partial M_{\leq^u}, \partial (M_{\leq^u})_t\big) < b\cdot t.
$$
\end{prop}

\begin{proof}
By definition $d(x,\partial M_{\leq^u})=t$ for any 
$x\in\partial (M_{\leq^u})_t$.

Let $x\in\partial \tilde M_{\leq^u}$ and let $c(t)$ be the geodesic line 
through $x$ with $\dot c(0)=\overrightarrow F(x)$ then
$$
 \frac{d}{dt}|_{t=0} D_{\{ d_\gc \leq\gep \}}\big(c(t)\big)=
 \overrightarrow F(x)\cdot \nabla D_{\{ d_\gc \leq\gep \}}(x)>1/b
$$ 
for any $\gc\in\Psi_{x,2\gep}$.
Since the convex function $D_{\{ d_\gc \leq\gep \}}\big( c(t)\big)$ has 
non-decreasing derivative we get that
$\frac{d}{dt}\Big( D_{\{ d_\alpha \leq\gep \}}\big( c(t)\big)\Big) >1/b$ 
for any $t>0$.
Thus, the point $c(bt)$ is outside  
$(\tilde M_{\leq^u})_t$, and so the distance from $x$ to
$\partial (\tilde M_{\leq^u})_t$ is less then $bt$. 
\end{proof}
 
The reason we prefer to work the $\gep$-neighborhood 
$M_{\leq^u}=(M_{\leq\gep})_\gep$ of the unipotent thin part is that it
has no ``sharp cusps''. 
In other wards, its curvature is uniformly bounded. 

\begin{lem}\label{smallC}
For an isometry $\gc$ and a point 
$x \in \partial (\{ d_{\gc} \leq \epsilon \} )_{\epsilon}$,
the $\gep$-ball for which the boundary sphere tangent at $x$ to the boundary of
$(\{ d_{\gc} \leq \epsilon \} )_{\epsilon}$, with the same external normal,
is contained in $(\{ d_{\gc} \leq \epsilon \} )_{\epsilon}$.
\end{lem}

\begin{proof}
The distance between $x$ and its closest point $\pi (x)$ in the closed convex 
set $\{d_\gc \leq \gep \}$ is easily seen to be $\gep$, and the $\gep$-ball 
centered at $\pi (x)$ is the required one.
(Recall that $\partial (\{d_\gc\leq\gep\} )_\gep$ is smooth.)
\end{proof}

The following two lemmas are clear.

\begin{lem}\label{B-B'}
Let $B'\subset B$ be topological spaces, and let $F_t~(t\in [0,1])$ be a 
deformation retract of $B$ such that $F_t(b)\in B'$ for any 
$b\in B',t\in [0,1]$. 
Then $F_t|_{B'}$ is a deformation retract of $B'$.
\end{lem}

Let $B_r$ denotes a ball of radius $r$ in $S$, and $B_r(x)$ the one centered
at $x\in S$.

\begin{lem}
There is a constant $l$ such that for any $\delta <1$, 
$$
 l\cdot \mu (B_{\delta /2})
 \geq \mu (B_{(b+1.5)\delta}).
$$ 
Thus any $\delta$-discrete subset of 
$(B_{(b+1)\delta})$ consists of at most $l$ elements.
\end{lem}

For a finite set
$\{ y_1,...,y_t\} \subset S$ we denote by $\sigma (y_1,...,y_t)$
its Chebyshev center, i.e. the unique point which minimize the function
$\max_{1\leq i\leq t}d(x,y_i)$.

We will soon take $B$ to be an intersection of balls, and $B'\subset B$ to be
the intersection of $B$ with the modified unipotent thick part 
$B'=B\cap M_{\geq^u}$. 
We intend to use \ref{B-B'} in order to show that, under some certain 
conditions $B'$ is contractible.
It is natural to use the so called star contraction to the Chebyshev center 
of the centers of the associated balls. This is the deformation retract which 
flows any point of $B$ on the geodesic segment which connect it to the 
required center, with constant velocity, 1 over the initial distance. 
In order to do this, we need the following:

\begin{prop}\label{positiveD}
There exist $0<\delta < \frac{\epsilon}{m(b+1)}$ 
such that for any point $x\in S$, any $\gep$-ball
$C$ which contains $x$ on its boundary sphere, and  
any $l$ points $y_1,\ldots ,y_l\in B_{(b+1)\delta }(x)\setminus 
(C)_{\delta}$,
the inner product of the external normal of $C$ at $x$ and the tangent at
$x$ to the geodesic segment $[x,\sigma (y_1,\ldots ,y_l)]$ is positive.
\end{prop}

\begin{proof}
The Riemannian metric on the ball of radius $(b+1)\gd$ is approximated by an 
ordinary Euclidean metric up to $o(\gd^2)$. Thus for small $\gd$, an 
approximated description of the above situation is given in the Euclidean 
space $\BR^{\dim (S)}$ by a half space, replacing the $\gep$-ball (as 
$\gep /\gd$ is very big), a point $x$ on its boundary hyper-plane and $l$ 
points at distance at least $\gd$ from this half space and at most $(b+1)\gd$ 
from $x$. Clearly, the unit vector pointing from $x$ to the Chebyshev center 
of these points has positive scalar product with the external normal to the 
half space at $x$, which is in fact $\geq$ then some positive constant 
depending only on $b$.

For more detailed proof see proposition 4.21 in \cite{Ge} and its proof.  
\end{proof}  

It follows from proposition \ref{closeN} that the union of a 
collection $\mathcal{C}$ of balls of radius 
$(b+1)\delta$, for which the set of centers form a maximal 
$\delta$-discrete subset of $M\setminus (M_{\leq^u} )_{\delta}$, covers 
$M_{\geq^u}=M\setminus M_{\leq^u}$. 
Note that $(b+1)\gd\leq\gep /m$, and thus, these balls are injected (see 
corollary \ref{InjRad}).
We fix such a collection, and for any subset $B\subset M$ we denote by
$B'$ its intersection with $M_{\geq^u}$.

Recall the following terminology from \cite{BT}: A cover of a topological space
$T$ is called {\it good cover} if any non-empty intersection of its sets is 
contractible. In that case, the simplicial complex which corresponds to the 
nerve of the cover (its vertices corresponds to the sets of the cover and a 
collection of vertices form a simplex if the corresponding sets contain 
a common point) is homotopically equivalent to $T$. 

\begin{prop}
Let $B$ be the intersection of $l$ (not necessarily deferent) balls of our 
collection $\mathcal{C}$, with centers $y_1,\ldots,y_l$. If $B$ is not empty 
then $\sigma (y_1,\ldots ,y_l)\in B'$ and the star-contraction from $B$ to
$\sigma (y_1,\ldots ,y_l)$ induces a contraction of $B'$. Hence the open 
set $B'$ is non-empty and diffeomorphic to $\BR^d$.
In particular  $\{ B':B\in \mathcal{C} \}$ is a good cover, and
the simplicial complex $\mathcal{R}$ corresponding to its nerve is 
homotopically equivalent to $M_{\geq^u}$, and therefore also to $M$.
\end{prop}

\begin{proof}
Lemma \ref{smallC} implies that for $\tilde x\in \partial\tilde M_{\leq^u}$ 
and for any  $\gc \in \gC^u$ for which 
$\tilde x\in \partial (\{ d_\gc \leq \epsilon \})_\epsilon$
there is an $\epsilon$-ball,
tangent to $\partial (\{ d_\gc \leq \epsilon \})_\epsilon$ at $\tilde x$,
which is contained in $(\{ d_\gc \leq \epsilon \})_\epsilon \subset 
\tilde M_{\leq^u}$.
Thus if in addition the image $x$ of $\tilde x$ belongs to $B$, proposition 
\ref{positiveD} implies that the geodesic segment 
$[x,\sigma (y_1,\ldots ,y_l)]$ is inside $B'$.
The proposition follows from lemma \ref{B-B'}
\end{proof}

We conclude that $M$ is homotopically equivalent to $\mathcal{R}$.
Since the collection of centers of the sets of $\mathcal{C}$ is $\gd$-discrete, 
we get that the number of vertices of $\mathcal{R}$, $|\mathcal{C}|$, is 
$\leq\frac{\mu (M)}{\mu (B_{\gd /2})}$. Since the sets of $\mathcal{C}$
are subsets of $(b+1)\gd$-balls, each of them intersects at most 
$\frac{\mu \big( B_{2(b+1.25)\gd}\big)}{\mu (B_{\gd /2})}$ of the 
others. Hence, all the vertices degrees of $\mathcal{R}$ are 
$\leq d:=\frac{\mu \big( B_{2(b+1.25)\gd}\big)}{\mu (B_{\gd /2})}$.
This complete the proof of theorem \ref{thm1}.
$\blacksquare$ 
\\

We now tern to the proof of \ref{thm2}. By Mostow's rigidity theorem, when
the universal covering space $S$ is not isometric to $\BH^2$, irreducible 
locally symmetric manifolds are characterized by their fundamental groups. 
If $M$ is non-compact arithmetic manifold locally isometric to $S$, and 
$\mathcal{R}$ is the simplicial complex corresponding to it by $\ref{thm1}$ 
then $\pi_1(M)=\pi_1(\mathcal{R})$.
The fundamental group of $\mathcal{R}$ is the same as that of its 2-skeleton
$\mathcal{R}^2$ (see \cite{Sp}). 
We use a combinatorial argument to estimate the number of possibilities for 
$\mathcal{R}^2$.
The 1-skeleton $\mathcal{R}^1$ is a graph with at most $c_1V$ vertices, where
$c_1=1/\mu (B_{\delta /2})$, and the degree of a vertices is at most 
$d=\mu (B_{2(b+1.25)\delta})/\mu (B_{\delta /2})$.

\begin{lem}
There are at most $V^{c_2V}$ such graphs.
\end{lem}

\begin{proof}
Going on the $\leq c_1V$ vertices one by one and choosing for each vertex a 
neighborhood from the available set of vertices at that stage, yield this 
estimation. 
\end{proof}

\begin{lem}
There are at most $c_1Vd^2$ triangles in such a graph.
\end{lem}

\begin{proof}
The number of triangles is bounded by the number of path of length 2, which
is at most $c_1Vd^2$. 
\end{proof}

\begin{cor}
There are at most $V^{CV}$ possibilities for the $2$-skeleton.
\end{cor}

\begin{proof}
Describing a 2-skeleton amounts to describing a 1-skeleton and choosing a 
subset of its set of triangles. 
Thus there are at most $V^{c_2V}2^{c_1Vd^2}\leq
V^{CV}$ possible 2-skeleton.
\end{proof}

This complete the proof of theorem \ref{thm2}.
$\blacksquare$

\end{document}